\documentclass[english]{amsart}

\usepackage{amsmath}
\usepackage{amssymb}
\usepackage{amsthm}
\usepackage{amscd}
\usepackage{babel}
\usepackage[latin1]{inputenc}
\usepackage{graphicx}

\renewcommand{\subjclassname}{AMS \textup{2000} Mathematics Subject
Classification\ }

\newtheorem{cor}{Corollary}
\newtheorem{prop}{Proposition}
\newtheorem{lem}{Lemma}
\theoremstyle{definition}

\author{Antonio M. Oller-Marc\'{e}n}
\title{A new look at the trailing zeroes of $N!$}
\address{Departamento de Matem\'{a}ticas, Universidad de Zaragoza\\
C/Pedro Cerbuna 12, 50009 Zaragoza (Espa\~{n}a)} \email{oller@unizar.es}
\begin{document}
\maketitle
\begin{abstract}
Let us denote by $Z_b(n)$ the number of trailing zeroes in the base $b$ expansi\'{o}n of $n!$. In this paper we study with some detail the behavior of the function $Z_b$. In particular, since $Z_b$ is non-decreasing, we will characterize the points where it increases and we will compute the amplitude of the jump in each of such points. In passing, we will study some asymptotic aspects and we will give families of integers that do not belong to the image of $Z_b$.
\end{abstract}
\subjclassname{11A25, 11A99}

\section{Introduction}
It is a usual exercise in Elementary Number Theory to compute the number of trailing zeroes in the base 10 expansion of the factorial of any integer. In fact, given any base $b$ and any integer $n$, it is easy to compute the number of trailing zeroes of the base $b$ expansion of $n!$. Namely, if we denote such number by $Z_b(n)$, we have.

\begin{lem}
\
\begin{enumerate}
\item $Z_p(n)=\displaystyle{\sum_{i\geq 1}\left[\frac{n}{p^i}\right]=\frac{n-\sigma_p(n)}{p-1}}$, where $\sigma_p(n)$ is the sum of the digits of the base $p$ expansion of $n$.
\item $\displaystyle{Z_{p^r}(n)=\left[\frac{Z_p(n)}{r}\right]}$ for every $r\geq 1$.
\item If $b=p_1^{r_1}\cdots p_s^{r_s}$, then $\displaystyle{Z_b(n)=\min_{1\leq i\leq s}Z_{p_i^{r_i}}(n)}$.
\end{enumerate}
\end{lem}

It is also easy to see that the function $Z_b:\mathbb{N}\longrightarrow\mathbb{N}$ is non-decreasing and not surjective. Thus, there exist integers $n$ such that $Z_b(n+1)>Z_b(n)$. In this situation we will say that $Z_b$ has a ``jump'' at $n$.

The previous lemma suggests the organization of the paper. The first section will be devoted to the prime case, the simplest one. We will characterize the points where $Z_p$ ``jumps'' and compute the amplitude of those ``jumps''. As an application we will give some families of integers not lying in the image of $Z_p$ and will study, in some sense, the density of $\textrm{Im}\ Z_p$. In the second section we will turn to the prime power case, where the results of the first section will be crucial. Lemma 1(3) shows that no further work is needed.

\section{The prime case}
During this section $p$ will be a prime. Given an integer $n$, let $n=\displaystyle{\sum_{i=0}^ka_ip^i}$ and $n+1=\displaystyle{\sum_{i=0}^ha'_ip^i}$ be the base $p$ expansions of $n$ and $n+1$ respectively. Also, let us define $t_{n,p}=\min\{j\ |\ a_j<p-1\}$. In other words, $a_0=\dots=a_{t_{n,p}-1}=p-1$ and $a_{t_{n,p}}<p-1$, being careful if $t_{n-p}=0$. If there is no risk of ambiguity we will just write $t$ instead of $t_{n,p}$.

In the following lemma, the relation between the digits of $n$ in base $p$ and tose of $n+1$ is studied. The proof is elementary and we omit.

\begin{lem}
$$a'_j=\begin{cases} 0 & \textrm{si $0\leq j< t$,}\\ a_t+1 & \textrm{si $j=t$,}\\ a'_j & \textrm{si $j>t$.}\end{cases}$$
\end{lem}

As a direct consequence we have the following proposition, which will be useful in the sequel.

\begin{prop}
\ 
\begin{enumerate}
\item $\sigma_p(n+1)-\sigma_p(n)=1-(p-1)t_{n,p}$.
\item $Z_p(n+1)-Z_p(n)=t_{n,p}$.
\end{enumerate}
\end{prop}
\begin{proof}
Part (1) follows directly from Lemma 2, while part (2) is a consequence of (1) together with Lemma 1(1).
\end{proof}

In the next result we characterize the points where $Z_p$ has a jump and we compute their amplitudes.

\begin{prop}
$Z_p(n+1)\neq Z_p(n)$ if and only if $p$ divides $n+1$. Moreover, $Z_p(n+1)-Z_p(n)=m$ if and only if $n+1=p^ma$ with $\textrm{g.c.d.}(p,a)=1$.
\end{prop}
\begin{proof}
It is enough to proof the second assertion. We will apply Proposition 1(2). Observe that $t=m$ if and only if $a_0=\dots =a_{m-1}=p-1>a_m$ which happens if and only if $a'_0=\dots=a'_{m-1}=0\neq a'_m$. This latter assertion is clearly equivalent to the fact that $p^m$ the greatest power of $p$ dividing $n+1$, and the proof is complete.
\end{proof}

Now, as an application of Proposition 2, we will give some families of integers not lying in $\textrm{Im}\ Z_p$. A first partial, but nevertheless interesting, result in this direction is the following.

\begin{cor}
The prime $p$ does not lie in $\textrm{Im}\ Z_p$; i.e., there is no $m\in\mathbb{N}$ such that the base $p$ expansion of $m!$ ends with $p$ zeroes.
\end{cor}
\begin{proof}
 By Lemma 1(1), $Z_p(p^2)=p+1$. Now, by Proposition 2, $Z_p(p^2-1)=p-1$ and the monotony of $Z_p$ completes the proof.
\end{proof}

In order to give a more general result in the style of the previous corollary, we will first need the following lemma which is a consequence of Lemma 1(1).
 
\begin{lem}
$\displaystyle{Z_p(lp^n)=l\frac{p^n-1}{p-1}+Z_p(l)}$.
\end{lem}

\begin{prop}
The following families of integers do not belong to $\textrm{Im}\ Z_p$:
\begin{itemize}
\item[a)] $\{\frac{p^n-kp+k-1}{p-1}\ |\ n>1,\ 1\leq k<n\}$.
\item[b)] $\left\{\left(\frac{p^k-1}{p-1}\right)p^n-k-h\ |\ n>1,\ k\geq1,\ 1\leq h<n\right\}$.
\end{itemize}
\end{prop}
\begin{proof}
Apply the previous Lemma with $l=1$ and $p^k-1$ respectively, together with the monotony of $Z_p$.
\end{proof}

Now we will study the density of $\textrm{Im}\ Z_p$. Given $N\in\mathbb{N}$, let us define $A_p(N)=\{n\leq N\ |\ n\in Z_p\}$ and put $a_p(N)=card(A_p(N))$. In the following proposition we study the asymptotic behavior of  $\frac{a_p(N)}{N}$.

\begin{prop}
$$\lim_{N\to\infty}\frac{a_p(N)}{N}=1.$$
\end{prop}
\begin{proof}
  Put $N=p^k-1$ and observe that $N\to\infty$ if and only if $k\to\infty$. We have that $Z_p((p-1)p^k)=N$ and, due to Proposition 2, until $(p-1)p^k$ they will take place $p-1$ ``jumps'' with amplitudes $1,\dots,k-1$. Consequently, $a_p(N)=N+1-\frac{(p-1)k(k-1)}{2}$ and
$$\frac{a_p(N)}{N}=1+\frac{1}{N}-\frac{(p-1)k(k-1)}{p^k-1}\longrightarrow 1$$ as desired.
\end{proof}

\section{The prime power case}
During this section $p$ will be a prime and $r\geq 2$ will be an integer. Recall that $\displaystyle{Z_{p^r}(n)=\left[\frac{Z_p(n)}{r}\right]}$, this fact will be crucial during this section.

\begin{prop}
 Let $n\in\mathbb{N}$ and put $Z_p(n)=\alpha r+\beta$. Then, $Z_{p^r}(n+1)=Z_{p^r}(n)$ if and only if $n+1=p^ma$ with $\textrm{g.c.d.}(p,a)=1$ and $0\leq m<r-\beta$.
\end{prop}
\begin{proof}
Let us suppose that $n+1=p^ma$ with $\textrm{g.c.d.}(p,a)=1$ and $0\leq m<r-\beta$. Then, by Proposition 2, we have that $Z_p(n+1)=Z_p(n)+m$. Thus, $Z_{p^r}(n+1)=\left[\frac{Z_p(n+1)}{r}\right]=\left[\frac{Z_p(n)+m}{r}\right]=\left[\frac{\alpha r+\beta+m}{r}\right]=\alpha=\left[\frac{Z_p(n)}{r}\right]=Z_{p^r}(n+1)$.

Conversely, assume that $Z_{p^r}(n+1)=Z_{p^r}(n)$. For some $m\geq 0$ it must be $Z_p(n+1)=Z_p(n)+m$ and thus,  $\alpha=\left[\frac{Z_p(n)}{r}\right]=\left[\frac{Z_p(n+1)}{r}\right]=\left[\frac{Z_p(n)+m}{r}\right]=\alpha+\left[\frac{\beta+m}{r}\right]$. From this, it follows that $\left[\frac{\beta+m}{r}\right]=0$. So $0\leq\beta+m<r$ and it is enough to recall Proposition 2 again to complete the proof.
\end{proof}

Now we will refine the previous result in order to compute the amplitude of the ``jumps''. It is interesting to compare the following result with Proposition 2.

\begin{prop}
Let $n\in\mathbb{N}$ and put $Z_p(n)=\alpha r+\beta$. Then, $Z_{p^r}(n+1)-Z_{p^r}(n)=k$ if and only if $n+1=p^ma$ with $\textrm{g.c.d.}(p,a)=1$ and $kr\leq m+\beta<(k+1)r$.
\end{prop}
\begin{proof}
Let us suppose that $Z_{p^r}(n+1)-Z_{p^r}(n)=k$. In such case, $Z_p(n+1)=Z_p(n)+m$ for some $m\geq 0$ and $\alpha+k=\left[\frac{Z_p(n)}{r}\right]+k=Z_{p^r}(n)+k=Z_{p^r}(n+1)=\left[\frac{Z_p(n+1)}{r}\right]=\left[\frac{Z_p(n)+m}{r}\right]=\left[\frac{\alpha r+\beta+m}{r}\right]=\alpha+\left[\frac{m+\beta}{r}\right]$.

Conversely, if $n+1=p^ma$ with $\textrm{g.c.d.}(p,a)=1$ and $kr\leq m+\beta<(k+1)r$, due to Proposition 2 we have that $Z_p(n+1)=Z_p(n)+m$, and consequently $Z_{p^r}(n+1)=\left[\frac{Z_p(n)}{r}\right]=\left[\frac{Z_p(n)+m}{r}\right]=\left[\frac{\alpha r+\beta+m}{r}\right]=\alpha+\left[\frac{m+\beta}{r}\right]=\alpha+k=\left[\frac{Z_p(n)}{r}\right]+k=Z_{p^r}(n)+k$.
\end{proof}

The rest of the section will be devoted to present families of integers not lying in $Z_{p^r}$ for various $p$ and $r$.

\begin{prop}
Let $p$ be a prime and $r\geq 2$ be an integer such that $r$ divides $\displaystyle{\sum_{i=1}^{kr}p^{kr-i}}$ for some $k>1$. Then, $\displaystyle{\frac{\sum_{i=1}^{kr}p^{kr-i}}{r}-h\notin\textrm{Im}\ Z_{p^r}}$ for every $1\leq h<k$.
\end{prop}
\begin{proof}
We have, $p$, $k$ and $r$ being like in the statement; that $Z_{p^r}(p^{kr})=\left[\frac{Z_p(p^{kr})}{r}\right]=\left[\frac{\sum_{i=1}^{kr}p^{kr-i}}{r}\right]=Z_{p^r}(p^{kr}-1)+k$. Again, the monotony of $Z_{p^r}$ completes the proof.
\end{proof}

The previous result can be slightly reformulated if either $r=2$ or $p=2$.

\begin{cor}
Let $p$ be an odd prime and let $k>1$ be an integer. Then, for every $1\leq h<k$,  $\displaystyle{\frac12\sum_{i=1}^{2k}p^{2k-i}}-h \notin Z_{p^2}$.
\end{cor}
\begin{proof}
It is enough to observe that $\displaystyle{\sum_{i=1}^{2k}p^{2k-i}}$ is even and apply the previous proposition.
\end{proof}

\begin{cor}
Let $q$ be an odd prime. Then, for every $1\leq h<q-1$,  $2^{q(q-1)}-1-h \notin Z_{2^q}$.
\end{cor}
\begin{proof}
By Fermat's little theorem $q$ divides $2^{q(q-1)}-1$, so we can take $k=q-1$ in Proposition 7 with $p=2$ and $r=q$.
\end{proof}

We will conclude the section and the paper with a new result in the style of the previous ones.

\begin{prop}
Let $p$ be a prime, $r\geq 2$ an integer and $l<p$ a multiple of $r$. Then, for every $k>1$ and $1\leq h<k$, $\displaystyle{\frac{l}{r}\sum_{i=1}^{kr}p^{kr-i}}-h\notin Z_{p^r}.$ 
\end{prop}
\begin{proof}
 It is clear that $Z_{p^r}(lp^{kr})=\frac{l}{r}\displaystyle{\sum_{i=1}^{kr}p^{kr-i}}$, while $Z_{p^r}(lp^{kr}-1)=Z_{p^r}(lp^{kr})-k$ so it is enough to apply again the monotony of $Z_{p^r}$.
\end{proof}

\end{document}